\documentclass{amsart}

\usepackage{latexsym}
\usepackage{amsmath, amssymb}
\usepackage[all]{xy}
\newtheorem{dfs}{Definition}[section]
\newtheorem{lms}[dfs]{Lemma}
\newtheorem{thms}[dfs]{Theorem}
\newtheorem{props}[dfs]{Proposition}

\newtheorem{cors}[dfs]{Corollary}
\newtheorem{rems}[dfs]{Remark}

\newcommand{\E}{\mathbf{Ell}}
\newcommand\VA{\mathrm{V}(A)}
\newcommand\WA{\mathrm{W}(A)}
\newcommand\W{\mathrm{W}}
\newcommand\V{\mathrm{V}}
\newcommand\TA{\mathrm{T}(A)}

\newcommand{\Z}{\ensuremath{\mathcal{Z}}}

\begin{document}

\title[Three Applications of the Cuntz Semigroup]{Three Applications of the Cuntz Semigroup}
\author{Nathanial P. Brown}
\address{Department of Mathematics, Penn State University, 
State College, PA, 16802, USA}
\email{nbrown@math.psu.edu}
\author{Andrew S. Toms}
\address{Department of Mathematics and Statistics, York University,
4700 Keele St., Toronto, Ontario, Canada, M3J 1P3}
\email{atoms@mathstat.yorku.ca}

\date{\today}

\thanks{N.B. was partially supported by DMS-0554870; A.T. was partially supported by NSERC}

\begin{abstract}
Building on work of Elliott and coworkers, we present three applications of the Cuntz semigroup:  
\begin{enumerate}
\item[(i)]  for many simple C$^*$-algebras, the Thomsen semigroup is recovered functorially from the Elliott 
invariant, and this yields a new proof of Elliott's classification theorem 
for simple, unital AI algebras; 
\item[(ii)] for the algebras in (i), classification of their Hilbert modules is similar to the von Neumann algebra context;    
\item[(iii)] for the algebras in (i), approximate unitary equivalence of self-adjoint operators is characterised
in terms of the Elliott invariant. 
\end{enumerate} 
\end{abstract}

\maketitle

\section{Introduction}

The Cuntz semigroup (cf.\ \cite{Cu}, \cite{KR}, \cite{pt}, \cite{Rfunct}  for definitions and basic properties) 
has recently become quite popular. In this note we extend the main theorem of \cite{bpt} to stable C$^*$-algebras. 
By combining this result with those of Coward-Elliott-Ivanescu \cite{cei} and Elliott-Ciuperca \cite{elliott-ciuperca}, 
we obtain the applications of the abstract directly. 
For the reader interested primarily in Elliott's classification program, we emphasize 
that most of our results are formulated in terms of the Elliott invariant -- the Cuntz 
semigroup is a powerful technical tool used only in proofs.    
This paper is a natural sequel to \cite{bpt} and \cite{pt}, and the latter contain
the requisite definitions, notation, and basic facts employed herein.

%%%%%%%%%%%%%%%%%%%%%%%%%%%%%%%%%%%%%%%%%%%%%%
\section{Computation of $\mathrm{W}(A\otimes \mathcal{K})$} 

Let $\WA$ denote the Cuntz semigroup of a unital simple C$^*$-algebra $A$ and let 
$\TA$ denote the simplex of tracial states on $A$.  We always assume that $\TA$ is 
nonempty, whence $A$ is stably finite.  It follows that $\WA$ can be decomposed 
into the disjoint union of $\VA$ (the Murray-von Neumann semigroup of equivalence classes
of projections) and the set $\WA_+$ of Cuntz classes of positive elements which are \emph{not} 
equal to the class of a projection. If $\mathrm{LAff}_b(\mathrm{T}(A))^{++}$ denotes the bounded, 
lower semicontinuous, affine, strictly positive functions on $\TA$, then there is a canonical map
\[
\iota\colon W(A)_+ \to \mathrm{LAff}_b(\mathrm{T}(A))^{++}
\]
given by
\[
\iota(\langle x\rangle)(\tau) = d_{\tau}(x), 
\]
where $d_{\tau}(x) := \lim_{n\to \infty} \tau\otimes \mathrm{Tr}_k (x^{1/n})$ for an element 
$x \in A\otimes M_k(\mathbb{C})$. (Here $\mathrm{Tr}_k$ is the \emph{non-normalized} trace on 
$M_k(\mathbb{C})$.) The main theorem of \cite{bpt} was that $\iota$ is an order isomorphism, whence 
\begin{equation}\label{structure1}
\WA \cong \mathrm{V}(A) \sqcup \mathrm{LAff}_b(\mathrm{T}(A))^{++}
\end{equation}
as partially ordered semigroups for two important classes 
of C$^*$-algebras: simple unital exact finite C$^*$-algebras
which absorb the Jiang-Su algebra $\mathcal{Z}$ tensorially, 
and simple unital AH algebras with slow dimension growth.  (We refer the reader 
to \cite{pt} for the definition of the order structure on $\mathrm{V}(A) \sqcup 
\mathrm{LAff}_b(\mathrm{T}(A))^{++}$.  As usual, $\mathcal{Z}$ denotes 
the Jiang-Su algebra -- see \cite{JS1}.) In this section we prove a structure theorem
similar to (\ref{structure1}) for $\mathrm{W}(A \otimes \mathcal{K})$, with $A$ as above.

Recall that $A$ has \emph{strict comparison} if $x \precsim y$  whenever 
$d_{\tau}(x) < d_{\tau}(y)$ for all $\tau \in \TA$. ($\precsim$ denotes 
Cuntz's relation and $x \sim y$ means $x \precsim y$ and $y \precsim x$.) 
When $A$ is unital simple exact and has strict comparison, the map $\iota$ 
is an isomorphism whenever it is surjective (cf.\ \cite[Proposition 3.3]{pt}). 

\begin{lms}\label{ctsapprox}
Let $A$ be a simple unital exact C$^*$-algebra with strict comparison of positive
elements and stable rank one.  Assume that $\iota$ is surjective (in particular,
the isomorphism of (\ref{structure1}) holds for $A$) and let $\langle a \rangle \in \WA_+$ be given. 
It follows that for any $\epsilon>0$, there exists $\delta>0$ and a continuous affine function
$f:\mathrm{T}(A) \to \mathbb{R}^+$ such that
\[
d_{\tau}((a-\epsilon)_+) < f(\tau) < d_{\tau}((a-\delta)_+), \ \forall \tau \in \mathrm{T}(A).
\]
\end{lms}

\begin{proof} Since $A$ has stable rank one, zero is an accumulation point of the spectrum $\sigma(a)$ (cf.\ \cite[Proposition 3.12]{perera}).  Choose points $\delta < \eta \in (0,\epsilon) \cap \sigma(a)$ so that
each of $(\delta,\eta)$ and $(\eta,\epsilon)$ are nonempty.
Since $A$ is simple, each trace and hence each lower semicontinuous dimension function is faithful.
It follows from a functional calculus argument that 
\[
d_{\tau}((a-\epsilon)_+) < d_{\tau}((a-\eta)_+) < d_{\tau}((a-\delta)_+), \ \forall \tau \in \mathrm{T}(A).
\]

Let $\mu_{\tau}$ be the (regular Borel) measure induced on $\sigma(a)$ by $\tau \in \mathrm{T}(A)$. 
The affine map $h:\mathrm{T}(A) \to \mathbb{R}^+$ given by
\[
h(\tau) := \mu_{\tau}( [\epsilon,\infty) \cap \sigma(a))
\]
is upper semicontinuous by the Portmanteau Theorem (\cite{Bi}).  From the inclusions
\[
(\epsilon,\infty) \cap \sigma(a) \subseteq [\epsilon, \infty) \cap \sigma(a) \subseteq (\eta,\infty) \cap \sigma(a)
\]
we have the following inequalities:
\[
d_{\tau}((a-\epsilon)_+) \leq h(\tau) \leq d_{\tau}((a-\eta)_+) < d_{\tau}((a-\delta)_+), \ \forall \tau \in \mathrm{T}(A).
\]
The affine map $\tau \mapsto d_{\tau}((a-\delta)_+)$ is lower semicontinuous.  Since $\mathrm{T}(A)$ is a metrizable
compact convex set, this map is the pointwise supremum of an increasing sequence of continuous affine maps, say $(f_n)_{n=1}^{\infty}$.
A straightforward argument using compactness then shows that there is some $n_0 \in \mathbb{N}$ such that
\[
f_n(\tau) > h(\tau), \ \forall \tau \in \mathrm{T}(A), \ \forall n \geq n_0.
\]
Setting $f(\tau) = f_{n_0}(\tau)$ completes the proof.
\end{proof}

Let $A$ be a unital C$^*$-algebra and $a \in A \otimes \mathcal{K}$
be positive.  Let $\{e_n\} \subset \mathcal{K}$ be an increasing sequence of 
projections with $\mathrm{rank}(e_n) = n$, and put $P_n = 1\otimes e_n \in A\otimes \mathcal{K}$. Then,
\[
P_1 a P_1 \precsim P_2 a P_2 \precsim P_3 a P_3 \precsim \cdots
\]
in $\W(A\otimes\mathcal{K})$ and $P_n a P_n \to a$ in norm.  Let $b = \sup_n \langle P_n a P_n \rangle \in \W(A\otimes\mathcal{K})$ 
(suprema of increasing sequences in the Cuntz semigroup always exist by \cite[Theorem 1]{cei}).
Then, given $\epsilon > 0$, there is some $n \in \mathbb{N}$ such that
\[
(a-\epsilon)_+ \precsim P_n a P_n \precsim b.
\]
It follows that $a \precsim b$ and $P_n a P_n \precsim a$ for each $n$.  Since the supremum is
unique, $a \sim b$.

\begin{lms}\label{specialsup} Let $A$  be as in Lemma \ref{ctsapprox}, and let 
$a \in A\otimes \mathcal{K}$ be a positive element such that $\langle a \rangle \in \W(A \otimes \mathcal{K})_+$. 
It follows that there is a sequence $(a_n)_{n=1}^{\infty}$ of positive elements in $A \otimes \mathcal{K}$
which satisfies the following conditions:
\begin{enumerate}
\item[(i)] $\langle a \rangle = \sup_n \langle a_n \rangle$;
\item[(ii)] $a_n \in A \otimes \mathrm{M}_{k(n)}$ for some $k(n) \in \mathbb{N}$;
\item[(iii)] the map $\tau \mapsto d_{\tau}(a_n)$ is continuous on $\mathrm{T}(A)$ for each $n \in \mathbb{N}$;
\item[(iv)] $d_{\tau}(a_n) < d_{\tau}(a_{n+1})$ for each $\tau \in \mathrm{T}(A)$ and $n \in \mathbb{N}$. 
\end{enumerate}
\end{lms}

\begin{proof}
Let $P_n$ be the unit of $A \otimes \mathrm{M}_n$ (as above) and define $b_n := P_n a P_n$.  
The sequence $b_n$ satisfies parts (i) and (ii) of the conclusion of the lemma
by construction.

\vspace{2mm}
\noindent
{\bf Case I.} Let us first address the case where infinitely many of the $b_n$s are Cuntz equivalent to 
a projection.  By passing to a subsequence, we may assume that every $b_n$ is  Cuntz equivalent 
to a projection (this does not affect the validity of (i) and (ii)).  If infinitely many of the
$b_n$s are Cuntz equivalent to a fixed projection $p \in A \otimes \mathcal{K}$, then we 
have 
\[
\langle a \rangle = \sup_n \langle b_n \rangle = \langle p \rangle;
\]
this contradicts our assumption that $a$ is not Cuntz equivalent to a projection.  Thus, 
each Cuntz class $\langle b_m \rangle$, $m \in \mathbb{N}$ occurs at most finitely many times
in the sequence $(\langle b_n \rangle)_{n=1}^{\infty}$.  Passing to a subsequence again,
we may assume that $\langle b_m \rangle \neq \langle b_n \rangle$ whenever $m \neq n$.  

Put $a_n = b_n$.  As noted, $(a_n)_{n=1}^{\infty}$ satisfies parts (i) and (ii) of the 
conclusion of the lemma already.  Condition (iii) is satisfied because the map $\tau 
\mapsto d_{\tau}(a_n)$ is continuous for any $a_n$ Cuntz equivalent to a projection.
Since $a_n$ is Cuntz equivalent to a projection, it is complemented inside 
$a_{n+1}$, i.e., there exists some positive element $c$ of $A \otimes \mathrm{M}_{k(n+1)}$

such that $\langle a_n \rangle + \langle c \rangle = \langle a_{n+1} \rangle$ (\cite[Proposition 2.2]{pt}).  Since $A$ is 
simple, $d_{\tau}(c) > 0$ for every $\tau \in \mathrm{T}(A)$.  It follows that $d_{\tau}(a_n) <
d_{\tau}(a_{n+1})$ for every $\tau$, as desired. 

\vspace{2mm}
\noindent
{\bf Case II.} Now we may assume that none of the $b_n$s are equivalent 
to a projection.  Given any $ \epsilon_n > 0$, we may use Lemma \ref{ctsapprox} to find $\delta_n>0$
and a continuous affine map $f_n:\mathrm{T}(A) \to \mathbb{R}^+$ such that
\[
d_{\tau}((b_n-\epsilon_n)_+) < f_n(\tau) < d_{\tau}((b_n-\delta_n)_+), \ \forall \tau \in \mathrm{T}(A).
\]
Choose inductively a sequence $(\epsilon_n)_{n=1}^{\infty}$ satisfying the following conditions:
\begin{enumerate}
\item[(a)] $(b_{n-1}-\delta_{n-1})_+ \precsim (b_n - \epsilon_n)_+$;
\item[(b)] $(b_k - \epsilon_k/n)_+ \precsim (b_n - \epsilon_n)_+$, $1 \leq k < n$.
\end{enumerate} 
By the surjectivity of $\iota$ there exists, for each $n \in \mathbb{N}$, a positive element $a_n$
in some $A \otimes \mathrm{M}_{k(n)}$ such that
\begin{equation}\label{dan}
d_{\tau}((b_n-\epsilon_n)_+) < f_n(\tau) = d_{\tau}(a_n) < d_{\tau}((b_n-\delta_n)_+), \ \forall \tau \in \mathrm{T}(A).
\end{equation}
The sequence $(a_n)_{n=1}^{\infty}$ therefore satisfies conditions (ii), (iii), and (iv) in the
conclusion of the lemma. 

Let us now verify condition (i).  By (b), (\ref{dan}), and the fact that $A$ has strict comparison of
positive elements we have $(b_k - \epsilon_k/n)_+ \precsim a_n$ for every $1 \leq k < n$ and $n \in \mathbb{N}$.
It follows that
\[
\sup_n \langle a_n \rangle \geq \langle (b_k-\epsilon_k/n)_+ \rangle, \ \forall n,k \in \mathbb{N}.
\]
In particular,
\[
\sup_n \langle a_n \rangle \geq \langle b_k \rangle, \ \forall k \in \mathbb{N},
\]
and so
\[
\sup_n \langle a_n \rangle \geq \sup_k \langle b_k \rangle = \langle a \rangle.
\]
On the other hand, $a_n \precsim b_n$ for every $n$, and so
\[
\sup_n \langle a_n \rangle \leq \sup_k \langle b_k \rangle = \langle a \rangle.
\]
This completes the proof.
\end{proof}

For each trace $\tau \in \TA$ and positive element $a \in A\otimes \mathcal{K}$ we define
a function $\iota\langle a\rangle\colon \mathrm{T}(A) \to \mathbb{R}^+ \cup \{\infty\}$ as follows:
\[
\iota\langle a\rangle(\tau) = \sup_n d_{\tau}(P_n a P_n).
\]
\begin{lms} 
\label{lms:independent}
$\iota\langle a\rangle$ is independent of the choice of projections $P_n$.  
\end{lms}  

\begin{proof}  Let $\{e_n\}, \{f_n\} \subset \mathcal{K}$ be increasing sequences of 
projections with $\mathrm{rank}(e_n) = \mathrm{rank}(f_n) = n$, and put $P_n = 1\otimes 
e_n, Q_n = 1\otimes f_n \in A\otimes \mathcal{K}$.  

Fix $n \in \mathbb{N}$ and $\epsilon > 0$.  Since $\lim_{k\to \infty} \|P_n Q_k - P_n\| =0$, 
we can find $k$ such that $\| P_n Q_kaQ_kP_n - P_naP_n \| < \epsilon$.  It follows that 
$(P_n a P_n - \epsilon)_+ \precsim Q_kaQ_k$ for all sufficiently large $k$.  In particular, 
$d_{\tau}((P_n a P_n - \epsilon)_+) \leq \sup_k d_{\tau}(Q_k a Q_k)$ for every $\epsilon > 0$.  
Since $d_{\tau}(P_naP_n) = \sup_{\epsilon}  d_{\tau}((P_n a P_n - \epsilon)_+)$, the lemma follows. 
\end{proof} 

\begin{props}\label{equivprop}
Let $A$ be a unital simple exact C$^*$-algebra with strict comparison and assume that 
$\iota\colon W(A)_+ \to \mathrm{LAff}_b(\mathrm{T}(A))^{++}$ is surjective.
Let $a,b \in A \otimes \mathcal{K}$ be positive elements which are not Cuntz equivalent to
projections.  It follows that $a \sim b$ if and only if $\iota\langle a\rangle = \iota\langle b\rangle$.
\end{props}

\begin{proof}
First suppose that $a \sim b$.  For each $\epsilon > 0$ and $n \in \mathbb{N}$ there exists 
a $\delta > 0$ and $m \in \mathbb{N}$ such that
\[
(P_n a P_n - 2 \epsilon)_+ \precsim (a-\epsilon)_+ \precsim (b-\delta)_+ \precsim P_m b P_m.
\]
It follows that for any $\tau \in \mathrm{T}(A)$,
\[
\iota\langle b\rangle(\tau) \geq d_{\tau}(P_m b P_m) \geq d_{\tau}(P_n a P_n - 2\epsilon)_+.
\]
Since $n$ and $\epsilon$ were arbitrary, we conclude that $\iota\langle b\rangle(\tau) \geq \iota\langle a\rangle(\tau)$.  Similarly, $\iota\langle a\rangle(\tau) \geq \iota\langle b\rangle(\tau)$.

Now suppose $\iota\langle a\rangle = \iota\langle b\rangle$.  Find, using Lemma \ref{specialsup}, sequences $(a_n)_{n=1}^{\infty}$ and
$(b_n)_{n=1}^{\infty}$ corresponding to $a$ and $b$, respectively.  By a compactness argument, for each $n \in \mathbb{N}$ there exists $m \in \mathbb{N}$ such that for every $\tau \in \mathrm{T}(A)$
we have the following inequalities:
\[
d_{\tau}(a_n) < d_{\tau}(b_m); \ \ d_{\tau}(b_n) < d_{\tau}(a_m).
\]
Since $A$ has strict comparison, $a_n \precsim b_m$ and 
$b_n \precsim a_m$.  It follows that
\[
\langle a \rangle = \sup_n \langle a_n \rangle = \sup_m \langle b_m \rangle = \langle b \rangle,
\]
as desired.  
\end{proof}

Let $\mathrm{SAff}(\mathrm{T}(A))$ denote the set of functions on $\mathrm{T}(A)$ which are pointwise
suprema of increasing sequences of continuous, affine, and strictly positive functions on $\mathrm{T}(A)$.  
Define an addition operation on the disjoint union $\VA \sqcup \mathrm{SAff}(\mathrm{T}(A))$ as follows:
\vspace{2mm}
\begin{enumerate}
\item[(i)] if $x,y \in V(A)$, then their sum is the usual sum in $V(A)$;
\item[(ii)] if $x,y \in \mathrm{SAff}(\mathrm{T}(A))$, then their sum
is the usual (pointwise) sum in $\mathrm{SAff}(\mathrm{T}(A))$;
\item[(iii)] if $x \in V(A)$ and $y \in \mathrm{SAff}(\mathrm{T}(A))$,
then their sum is the usual (pointwise) sum of $\hat{x}$ and $y$ in
$\mathrm{SAff}(\mathrm{T}(A))$, where $\hat{x}(\tau) = \tau(x)$,
$\forall \tau \in \mathrm{T}(A)$.
\end{enumerate}
\vspace{2mm}
Equip $\VA \sqcup \mathrm{SAff}(\mathrm{T}(A))$ with the partial order $\leq$ which restricts to the
usual partial order on each of $V(A)$ and $\mathrm{SAff}(\mathrm{T}(A))$,
and which satisfies the following conditions for $x \in V(A)$ and $y \in
\mathrm{SAff}(\mathrm{T}(A))$:
\vspace{2mm}
\begin{enumerate}
\item[(i)] $x \leq y$ if and only if $\hat{x}(\tau) < y(\tau)$, $\forall \tau \in \mathrm{T}(A)$;
\item[(ii)] $y \leq x$ if and only if $y(\tau) \leq \hat{x}(\tau)$, $\forall \tau \in \mathrm{T}(A)$.
\end{enumerate}
\vspace{2mm}

\begin{thms}\label{custructure}
Let $A$ be a unital simple exact and tracial C$^*$-algebra with strict comparison. 
Assume that $\iota\colon W(A)_+ \to \mathrm{LAff}_b(\mathrm{T}(A))^{++}$ is surjective. 
It follows that $$\W(A\otimes\mathcal{K}) \cong \VA \sqcup \mathrm{SAff}(\mathrm{T}(A)),$$ as ordered semigroups. 
\end{thms}

\begin{proof}  
Define 
\[
\phi:\W(A \otimes \mathcal{K}) \to \VA \sqcup \mathrm{SAff}(\mathrm{T}(A))
\]
by $\mathbf{id}_{\V(A \otimes \mathcal{K})}$ on $\V(A \otimes \mathcal{K})$
and by $\iota$ on $\W(A \otimes \mathcal{K})_+$.  Let us first prove that $\phi$
is a bijection.
For each $x \in \W(A\otimes\mathcal{K})_+$ we define an element 
$\iota (x) \in \mathrm{SAff}(\mathrm{T}(A))$ as follows: choose a positive element 
$a \in M_n(A\otimes\mathcal{K}) \cong A\otimes\mathcal{K}$ such that $x = \langle a \rangle$ 
and define $\iota (x) := \iota \langle a \rangle$ (as in the previous proposition).  
(Since $\iota \langle a \rangle$ is independent of the projections used in its 
definition, it is not hard to check that the definition of $\iota (x)$ is independent 
of the identification $M_n(\mathcal{K}) \cong \mathcal{K}$.) Since $A$ is stably finite, 
it suffices to prove that $\iota\colon \W(A\otimes\mathcal{K})_+ \to \mathrm{SAff}(\mathrm{T}(A))$ 
is a bijection. 

Injectivity of $\iota$ follows from Proposition \ref{equivprop}.  Surjectivity follows from 
two facts: (1) the range of $\iota$ contains $\mathrm{LAff}_b(\mathrm{T}(A))^{++}$ and (2) 
$\W(A\otimes\mathcal{K})$ has suprema (cf.\ \cite{cei}). Indeed, given 
$f \in \mathrm{SAff}(\mathrm{T}(A))$ we find continuous affine functions $f_n \leq f_{n+1}$ 
converging up to $f$ pointwise.  Letting $a_n \in A\otimes\mathcal{K}$ be positive elements 
such that $\hat{a}_n = f_n$, we let $x = \sup_n \langle a_n \rangle \in \W(A\otimes\mathcal{K})$ 
(we have used strict comparison here to ensure $\{\langle a_n \rangle\}$ is an increasing 
sequence in $\W(A\otimes\mathcal{K})$). Then it is clear that $\iota(x) = f$.

To complete the proof, we must show that $\phi$ is order preserving.  Suppose that 
$x \leq y$, $x,y \in \W(A \otimes \mathcal{K})$.  There are four cases to consider.
\begin{enumerate}
\item[(a)] If $x,y \in \V(A \otimes \mathcal{K})$,
then $\phi(x) \leq \phi(y)$ since $\phi|_{\V(A \otimes \mathcal{K})} = \mathbf{id}_{\V(A \otimes \mathcal{K})}$.
\item[(b)] If $x,y \in \W(A \otimes \mathcal{K})$, then $\phi(x) \leq \phi(y)$ since
$\phi|_{\W(A \otimes \mathcal{K})} = \iota$ and $\iota$ is order-preserving.
(The proof of this last fact follows from the proof of the first implication
in Proposition \ref{equivprop}.)  
\item[(c)] If $x \in \V(A \otimes \mathcal{K})$ and $y
\in \W(A \otimes \mathcal{K})_+$, then we apply \cite[Proposition 2.2]{pt} to
find $z \in \W(A \otimes \mathcal{K})$ such that $x + z = y$.  It follows that
$\iota(x)(\tau) < \iota(y)(\tau)$, $\forall \tau \in \mathrm{T}(A)$ (note that
$\iota(x)(\tau) < \infty$ in this case), whence $\phi(x) \leq \phi(y)$. 
\item[(d)] If $x \in \W(A \otimes \mathcal{K})_+$ and $y \in \V(A \otimes \mathcal{K})$, 
then $\phi(x) \leq \phi(y)$ since $\iota$ is order-preserving.
\end{enumerate}
\end{proof}

The theorem above holds for all simple unital 
AH algebras with slow dimension growth, and for the class of simple unital 
exact stably finite C$^*$-algebras which absorb $\mathcal{Z}$ (\cite{bpt}, \cite{T1}). 

%%%%%%%%%%%%%%%%%%%%%%%%%%%%%%%%%%%%%%%%%%%%%%
\section{Classifying Hilbert Modules}

Let $E,F $ be countably generated Hilbert modules over a separable, unital C$^*$-algebra $A$.  
By Kasparov's stabilization theorem, there are projections $P_E, P_F \in L(H_A)$ such that 
$E$ is isomorphic to $P_E H_A$ and $F$ is isomorphic to $P_F H_A$. (Here $H_A = \ell^2 \otimes A$ 
is the standard Hilbert module over $A$ and $L(H_A)$ is the set of bounded adjointable operators on 
$H_A$; see \cite{lance} for more).  Since $L(H_A) = M(A\otimes \mathcal{K})$ (the multiplier algebra 
of $A\otimes \mathcal{K}$), we can find strictly positive elements $a \in P_E(A\otimes \mathcal{K})P_E$ 
and $b \in P_F(A\otimes \mathcal{K})P_F$. According to \cite[Theorem 3]{cei}, if we further assume $A$ 
has stable rank one, $$E\cong F \mbox{ if and only if } \langle a \rangle = \langle 
b \rangle \in \W(A\otimes \mathcal{K}).$$ 
In this section we'll reformulate this result in terms of the projections $P_E$ and $P_F$. 

First, an alternate formula for $\iota\langle a \rangle \in \mathrm{SAff}(\mathrm{T}(A))$ will be handy. 
Let $\mathcal{F} \subset \mathcal{K}$ denote the \emph{finite-rank operators} and $A\otimes \mathcal{F}$ 
be the algebraic tensor product of $A$ and $\mathcal{F}$ (which we identify with the ``finite-rank" 
operators on $H_A$).

\begin{lms}  For every $0 \leq a \in A\otimes \mathcal{K}$ and $\tau \in \TA$ we have
$$\iota\langle a\rangle(\tau) = \sup \{d_{\tau}(b) : 0 \leq b \in A\otimes \mathcal{F}, b \precsim a\}.$$ 
\end{lms} 

\begin{proof}  If $P = 1 \otimes e$ for some finite rank projection $e \in \mathcal{K}$, then $P a P  \in A\otimes
 \mathcal{F}$ and $P a P \precsim a$; hence, the inequality $\leq$ is immediate.  

For the other direction, fix $b \in A\otimes \mathcal{F}$ such that $b \precsim a$, and fix $\epsilon > 0$.  
Choose $\delta > 0$ such that $d_{\tau}(b) - \epsilon \leq d_{\tau}((b-\delta)_+)$ and find $x \in A\otimes 
\mathcal{K}$ such that $\| x^* a x - b \| < \delta$. By density, we may assume $x \in A\otimes M_n(\mathbb{C})$ 
for some large $n \in \mathbb{N}$. It follows that $(b-\delta)_+ \precsim x^* a x$. Now, let $P_n = 1\otimes e_n$, 
for some increasing finite-rank projections $e_n$, such $P_n x = x = xP_n$ for all $n$.  We have that 
$(b-\delta)_+ \precsim x^* a x = x^*(P_n aP_n)x$.  Hence, $$d_{\tau}(b) - \epsilon \leq d_{\tau}((b-\delta)_+) 
\leq d_{\tau}(P_n a P_n),$$ and, by Lemma \ref{lms:independent}, this completes the proof of the lemma. 
\end{proof} 

For any projection $Q \in M(A\otimes \mathcal{K})$ and tracial state $\tau \in \TA$ we define $$\hat{Q}(\tau) = 
\sup \{ \tau \otimes \mathrm{Tr}(b) : 0 \leq b \in A\otimes \mathcal{F}, b \leq P\},$$ where $\mathrm{Tr}$ is the 
(unbounded) trace on $\mathcal{F}$. 

\begin{lms}  Assume $A$ is unital with stable rank one. For any projection $Q \in M(A\otimes \mathcal{K})$, 
strictly positive element $a \in Q(A\otimes \mathcal{K})Q$ and $\tau \in \TA$, we have $$\hat{Q}(\tau) = 
\iota\langle a\rangle(\tau).$$
\end{lms} 

\begin{proof}  Since $\{b: b \in A\otimes \mathcal{F}, b \leq P\} \subset \{b: b \in A\otimes \mathcal{F}, b \precsim a\}$ 
(cf.\ \cite[Propostion 2.7(ii)]{KR}), and $\tau \otimes \mathrm{Tr}(b) \leq \lim_n \tau \otimes \mathrm{Tr}(b^{1/n}) = 
d_{\tau}(b)$, the previous lemma implies that  $\hat{Q}(\tau) \leq \iota\langle a\rangle(\tau).$

For the opposite inequality, fix $b \in A\otimes \mathcal{F}$ such that $b \precsim a$, and $\epsilon > 0$.  
Choose $\delta > 0$ such that $d_{\tau}(b) - \epsilon \leq d_{\tau}((b-\delta)_+)$.  Since $A$ has stable rank 
one, so does $(A\otimes \mathcal{K}\tilde{)}$ (the unitzation of $A\otimes \mathcal{K}$). Hence, by 
\cite[Proposition 2.4]{Rfunct}, we can find a unitary $u \in (A\otimes \mathcal{K}\tilde{)}$ such that 
$u^* (b-\delta)_+ u \leq Q$.  Since $u^* (b-\delta)_+ u \in A\otimes \mathcal{F}$, the following inequalities 
complete the proof: $$d_{\tau}(b) - \epsilon \leq d_{\tau}((b-\delta)_+) = d_{\tau}(u^* (b-\delta)_+ u)  = 
\lim_n \tau \otimes \mathrm{Tr}([u^* (b-\delta)_+ u]^{1/n}) \leq \hat{Q}(\tau).$$
\end{proof} 

Recall that if $M \subset B(L^2(M))$ is a II$_1$-factor in standard form, then isomorphism classes of modules over $M$ (i.e.\ normal representations $M\subset B(H)$) are completely determined by the traces of the corresponding projections in $M^\prime \overline{\otimes} B(H)$.  Our next theorem is analogous to this classical result. 

\begin{thms}\label{moduleclass} Let $A$ be a unital simple exact and tracial C$^*$-algebra with strict comparison and stable rank one.  Assume that $\iota\colon W(A)_+ \to \mathrm{LAff}_b(\mathrm{T}(A))^{++}$ is surjective.  Given two countably generated Hilbert modules $E$, $F$ over $A$, the following are equivalent: 
\begin{enumerate} 
\item $E$ is isomorphic to $F$; 

\item $P_E$ is Murray-von Neumann equivalent to $P_F$; 

\item Either $\langle P_E\rangle = \langle P_F\rangle \in \mathrm{V}(A)$ (in the case 
$P_E, P_F \in A\otimes \mathcal{K}$), or $\hat{P}_E = \hat{P}_F$. 
\end{enumerate} 
In particular, if neither $E$ nor $F$ is a finitely generated projective module, then $E \cong F$ 
if and only if $\hat{P}_E = \hat{P}_F$.   
\end{thms} 

\begin{proof}  The equivalence of the first two conditions is well-known (and a simple exercise).  

Since we already mentioned that $E\cong F$ if and only if $\langle a \rangle = \langle 
b \rangle \in \W(A\otimes \mathcal{K})$ (\cite[Theorem 3]{cei}), where $a$ (resp.\ $b$) is a strictly 
positive element in $P_E(A\otimes \mathcal{K})P_E$ (resp.\ $P_F (A\otimes \mathcal{K})P_F$), the 
equivalence of (1) and (3) follows from the previous lemma and Theorem \ref{custructure}.
\end{proof} 

\begin{rems} {\rm The theorem above is, in a certain sense, best possible: we really need strict comparison.  More precisely, the hypotheses are satisfied by simple AH algebras with slow dimension growth (and $\mathcal{Z}$-stable algebras --  cf.\ \cite{bpt} and Theorem \ref{custructure}), but the result can't be extended to all AH algebras.  Indeed, the reader will find in \cite{T2} a pair of positive elements in a simple unital AH algebra of stable rank one  such that the corresponding Hilbert modules, say $E$ and $F$, are not isomorphic but do satisfy $\hat{P}_E = \hat{P}_F$.   

It is also worth remarking that the result above gives a complete parametrization of isomorphism classes of countably generated Hilbert modules over $A$ in terms of $\mathrm{K}_0$ and traces.}
\end{rems}

%%%%%%%%%%%%%%%%%%%%%%%%%%%%%%%%%%%%%%%%%%%%%%
\section{From Elliott to Thomsen and the Classification of Simple AI Algebras}

\begin{thms}\label{ell-thom} Let $A$ be a unital simple C$^*$-algebra of stable rank one 
for which $\mathrm{W}(A\otimes \mathcal{K}) = \VA \sqcup \mathrm{SAff}(\mathrm{T}(A))$. 
Then, the Thomsen semigroup of $A$ can be functorially recovered from the Elliott invariant of $A$. 
\end{thms} 

This theorem follows immediately from \cite{elliott-ciuperca}.  The result applies
to any algebra satisfying the hypotheses of Theorem \ref{custructure} -- in particular, 
by \cite{bpt}, $A$ could be a  simple unital AH algebra with slow dimension growth, or 
a simple unital exact and $\mathcal{Z}$-stable C$^*$-algebra.  The assumption of simplicity
in the theorem is actually redundant.  The assumption on the structure
of $\W(A \otimes \mathcal{K})$ gurantees that every trace on $A$ is faithful, whence $A$ is simple.

\begin{thms}[Elliott, \cite{elliott-AI}] Let $A$ and $B$ be simple unital 
inductive limits of algebras of the form $F\otimes C[0,1]$, where $F$ is 
finite dimensional.  Then $A\cong B$ if and only if $\E(A) \cong \E(B)$. 
\end{thms} 

\begin{proof} If $\E(A) \cong \E(B)$ then $\mathrm{W}(A\otimes \mathcal{K}) 
\cong \mathrm{W}(B\otimes \mathcal{K})$, by Theorem \ref{custructure} and \cite{bpt} 
(since AI algebras have no dimension growth). From \cite{elliott-ciuperca} it follows 
that the Thomsen semigroups of $A$ and $B$ are isomorphic too.  Hence, 
by \cite[Theorem 1.5]{Th}, $A \cong B$. 
\end{proof} 

Evidently, the result can be slightly improved -- the presence of $K_1$ is irrelevant. 
This theorem is best possible in the sense that the Elliott invariant 
is not complete for \emph{non-simple} AI algebras (cf.\ \cite{Th}). 
The Cuntz semigroup, however, is a complete invariant in the non-simple case,
as shown in \cite{elliott-ciuperca}.

\section{Unitary Orbits of Self-Adjoints in Simple, Unital, Exact C$^*$-algebras}

Let $a \in A$ be self-adjoint with spectrum $\sigma(a)$. 
Let $\phi_a:\mathrm{C}(\sigma(a)) \to A$ be the canonical homomorphism induced by sending the 
generator $z$ of $\mathrm{C}(\sigma(a))$ to $a \in A$, and denote by $\E(a)$ the following 
pair of induced maps:
\[
\mathrm{K}_*(\phi_a):\mathrm{K}_*(C(\sigma(a))) \to \mathrm{K}_*(A); \ \  \phi_a^{\sharp}: \mathrm{T}(A) \to \mathrm{T}(C(\sigma(a))).
\]
As in Theorem \ref{ell-thom}, the hypotheses of the next theorem guarantee the
simplicity of $A$.

\begin{thms}\label{unit-equiv}  Let $A$ be a unital C$^*$-algebra of stable rank one, and assume that 
\[
\W(A\otimes\mathcal{K}) \cong \VA \sqcup \mathrm{SAff}(\mathrm{T}(A)). 
\]
Let $a, b \in A$ be self-adjoint.  Then, $a$ and $b$ are approximately unitarily equivalent if and only 
if $\sigma(a) = \sigma(b)$ and $\E(a) = \E(b)$. 
\end{thms} 

\begin{proof}  The ``only if" statement is routine, so assume $\sigma(a) = \sigma(b)$ and $\E(a) = \E(b)$.  

First, we handle the case that $\sigma(a) = \sigma(b) \subset (0,\infty)$, i.e., that both $a$ and $b$ 
are positive and invertible.  Let $X = \sigma(a) = \sigma(b)$ and $\W_a\colon \W(C(X)) \to \W(A\otimes 
\mathcal{K})$ (resp.\ $\W_b\colon \W(C(X)) \to \W(A\otimes \mathcal{K})$) denote the Cuntz-semigroup map 
induced by the canonical homomorphism $C(X) \to A\otimes \mathcal{K}$ sending $z \mapsto a\otimes e_{1,1}$ 
(resp.\ $z \mapsto b\otimes e_{1,1}$).  We claim that $\W(a) = \W(b)$. 

So, let $h \in M_n(C(X))$ be positive and $h_a \in M_n(A)$ (resp.\  $h_b \in M_n(A)$) denote the image of 
$h$ under the canonical inclusion $M_n(C(X)) \subset M_n(A)$ sending $C(X) \to C^*(a)$ (resp.\ $C(X) \to 
C^*(b)$). Since $A$ and $C(X)$ have stable rank one, \cite[Proposition 3.12]{perera} implies that $h$ is 
equivalent to a projection in $M_{\infty}(C(X))$ if and only if $h_a$ is equivalent to a projection in 
$M_{\infty}(A)$. Thus, in this case, $\W_a\langle h \rangle = \langle h_a \rangle = \langle h_b \rangle = 
\W_b\langle h \rangle$ since $\langle h_a \rangle$ and $\langle h_b \rangle$ are the same element in 
$\mathrm{V}(A\otimes \mathcal{K}) \subset \W(A\otimes \mathcal{K})$. 

If $h$ is not equivalent to a projection -- i.e.\ if $\langle h_a \rangle, \langle h_b \rangle \in 
\mathrm{SAff}(\mathrm{T}(A)) \subset \W(A\otimes\mathcal{K})$ -- then it suffices to show that $d_{\tau}(h_a) 
= d_{\tau}(h_b)$ for every $\tau \in \TA$. However, if $\mu$ is a measure on $\sigma(h)$ then $d_{\mu}(h) 
= \mu(\sigma(h)\setminus \{0\})$.  Since $\E(a) = \E(b)$, the maps on tracial spaces agree -- i.e.\ for 
each $\tau \in \TA$ the measures induced by restriction agree on $\sigma(h_a) = \sigma(h_b)$ -- and hence 
$d_{\tau}(h_a) = d_{\tau}(h_b)$ for every $\tau \in \TA$, as desired. 

Knowing that $\W_a = \W_b$, it now follows from \cite{elliott-ciuperca} that $a\otimes e_{1,1}$ is 
approximately unitarily equivalent to $b\otimes e_{1,1}$ in the unitization of $A\otimes \mathcal{K}$. 
So, let $v_n \in (A\otimes \mathcal{K})^+$ be unitaries such that $v_n(a\otimes e_{1,1})v_n^* \to b\otimes e_{1,1}$.  
Since $a$ is invertible, for every $\varepsilon > 0$ there exists a polynomial $p$ such that $\|p(a) - 1\| 
< \varepsilon$; since $\sigma(a) = \sigma(b)$, $\|p(b) - 1\| < \varepsilon$ as well.  Hence, for large $n$, 
$\|v_n(1\otimes e_{1,1})v_n^* - 1\otimes e_{1,1}\| < C\varepsilon$ for some constant $C$ depending only on 
$\sigma(a)$.  If $\varepsilon$ is sufficiently small, this implies that $(1\otimes e_{1,1})v_n(1\otimes e_{1,1})$ 
is almost a unitary in $A$ -- hence can be perturbed to an honest unitary $u_n$.  A routine exercise now confirms 
that $a$ is approximately unitarily equivalent to $b$ (in $A$). 

For the case of general self-adjoints $a,b \in A$, we deduce the theorem from a simple trick.  Namely, 
fix some constant $c$ such that $a + c1$ is positive and invertible.  Then $b + c1$ is also positive 
and invertible.  By the case handled above, $a + c1$ and $b + c1$ are approximately unitarily equivalent, 
hence the same is true of $a$ and $b$. 
\end{proof}

The theorem above holds for all simple unital 
AH algebras with slow dimension growth, and for the class of simple unital 
exact stably finite $\mathcal{Z}$-stable C$^*$-algebras (see \cite{bpt}, \cite{T1}, 
and Theorem \ref{custructure}). 

A more general version of Theorem \ref{unit-equiv} holds for simple unital exact 
and stably finite C$^*$-algebras
\begin{cors} 
Let $a$ and $b$ be self-adjoint elements of a simple unital exact and stably finite
C$^*$-algebra $A$.  Then $a$ and $b$ are approximately unitarily equivalent in $A\otimes \mathcal{Z}$ --  
there exist unitaries $u_n \in A\otimes \mathcal{Z}$ such that $\| u_n (a\otimes 1) u_n^* - b\otimes 1 \| \to 0$ 
-- if and only if $\sigma(a) = \sigma(b)$ and $\E(a) = \E(b)$. 
\end{cors}

\noindent
The proof of this corollary is a tiny perturbation of the proof of Theorem \ref{unit-equiv}. The result is also, in some sense, best-possible:  in \cite{T2} a pair of positive elements in a simple unital
AH algebra were constructed which have identical Elliott data but which are not Cuntz equivalent 
(hence not unitarily equivalent).  For the truly interested reader, the elements in question are
$f(\tau^*(\xi) \times \tau^*(\xi))$ and $f \theta_1 \oplus f \theta_1$, constructed in Section 3
of \cite{T2}.


\begin{thebibliography}{999}

\bibitem{Bi} Billingsley, P.: {\it Convergence of Probability measures}, J. Wiley \& Sons, 1968.

\bibitem{bpt} Brown, N.P., Perera, F., and Toms, A. S.: {\it The Cuntz semigroup, the Elliott conjecture, and dimension
functions on C$^*$-algebras}, arXiv preprint math.OA/0609182 (2006) 

\bibitem{elliott-ciuperca} Ciuperca, A., and Elliott, G. A.: {\it A remark on invariants for C*-algebras of stable rank one}, in preparation. 

\bibitem{cei} Coward, K.T., Elliott, G. A., and Ivanescu, C.: \emph{The Cuntz semigroup as an invariant for C$^*$-algebras}, preprint (2006). 

\bibitem{Cu} Cuntz, J.: \emph{Dimension functions on simple C$^*$-algebras}, Math. Ann., \textbf{233} (1978), pp. 145--153.

\bibitem{elliott-AI} Elliott, G.A.: \emph{A classification of certain simple $C\sp *$-algebras},  
Quantum and non-commutative analysis (Kyoto, 1992),  373--385, Math.\ Phys.\ Stud., 16, Kluwer Acad. Publ., Dordrecht, 1993.

\bibitem{JS1} Jiang, X.\ and Su, H.: {\it On a simple unital projectionless $C^{*}$-algebra}, Amer.\ J.\ Math.\ {\bf 121} (1999), pp. 359-413.

\bibitem{KR} Kirchberg, E.\ R\o rdam, M.: \emph{Non-simple purely infinite C$^*$-algebras}, Amer. J. Math., \textbf{122} (2000), pp.
637--666.

\bibitem{lance} Lance, E.C.: \emph{Hilbert $C\sp *$-modules. A toolkit for operator algebraists}, 
London Mathematical Society Lecture Note Series, 210. Cambridge University Press, Cambridge, 1995. x+130 pp.

\bibitem{perera} Perera, F.: \emph{The structure of positive elements for C$^*$-algebras with real rank zero}, 
Inter.\ J.\ Math.\ \textbf{8} (1997), pp. 383--405. 

\bibitem{pt} Perera, F.\ and Toms, A. S.: \emph{Recasting the Elliott Conjecture}, Math. Ann.,
published online at http://dx.doi.org/10.1007/s00208-007-0093-3, pp. 1-34.

\bibitem{Rfunct} R\o rdam, M.: \emph{On the structure of simple C$^*$-algebras tensored with a UHF-algebra. II}, 
J. Funct. Anal., \textbf{107} (1992), pp. 255--269.

\bibitem{R1} R\o rdam, M.: \emph{The stable and the real rank of \Z-absorbing C$^*$-algebras}, International 
J. Math., \textbf{15} (2004), pp. 1065--1084.


\bibitem{Th} Thomsen, K.: {\it Inductive limits of interval algebras:  unitary orbits of positive elements},  
Math. Ann. {\bf 293} (1992), pp. 47-63.

\bibitem{T1} Toms, A. S.: {\it Stability in the Cuntz semigroup of a commutative C$^*$-algebra}, Proc. London Math. Soc., to appear,
arXiv preprint math.OA/0607099 (2006).

\bibitem{T2} Toms, A. S.: {\it On the classification problem for nuclear C$^*$-algebras}, Ann. of Math. (2), to appear,
arXiv preprint math.OA/0509103 (2005).

\end{thebibliography}
\end{document}